\documentclass[12pt,a4paper]{article}
\usepackage{amsmath, amssymb, amsthm}
\usepackage{mathtools}
\usepackage{microtype}
\usepackage[utf8]{inputenc}
\usepackage[T1]{fontenc}
\usepackage{geometry}
\usepackage{times}
\usepackage{natbib}
\usepackage{hyperref}
\hypersetup{
    colorlinks=true,
    linkcolor=blue,
    filecolor=magenta,      
    urlcolor=blue,
    citecolor=blue,
}

\newtheorem{theorem}{Theorem}[section]
\newtheorem{corollary}[theorem]{Corollary}

\title{Beyond Wald's Equation and the\\Optional Sampling Theorem}
\author{Michael J. Klass and Victor H. de la Pe\~{n}a}
\date{July 1, 2026}

\begin{document}

\maketitle

\begin{abstract}
This paper establishes a conservation identity for mean-zero martingales stopped by extended-valued stopping times. For any mean-zero martingale $\{M_n\}$ and any extended-valued stopping time $T$ satisfying $E|M_T|I(T<\infty)<\infty$, the quantity $L\equiv E[M_T I(T<\infty)]$ exists and equals $\lim_n E[-M_n I(T>n)]$, a limit which always exists. The optional sampling theorem for stopping times and uniformly integrable martingales---and Wald's equation for mean-zero random variables, as its i.i.d.\ specialization---is recovered with a little extra effort, in which case the limit also vanishes. The identity itself remains in force whether or not $L=0$, and whether or not $P(T<\infty)=1$. Two corollaries and an application derived from this identity provide information on the rate of decay of the tail probability of the stopping time.
Moreover, a necessary and sufficient condition is presented to characterize when $E|M_T|I(T<\infty)$ is finite. The characterization applies more generally whenever $|M_n|$ is a sequence of random variables, each having finite expectation.

A third theorem provides sufficient conditions ensuring that certain exceedance-level, potentially extended-valued, stopping times are finite with probability one. It further implies that $\limsup M_n=\infty$ almost surely. We demonstrate these results through examples and explore their implications for different families of martingales. Our findings extend classical results in martingale theory and provide new insights into the behavior of stopped martingales, especially when the expected value of the stopped martingale on the set where the extended-valued stopping time $T$ is finite differs from the expected value of the martingale at time 1.

\bigskip\noindent
{\bf Keywords}: Martingale; stopping time; Wald's equation; optional sampling theorem; randomly stopped sum.

\smallskip\noindent
{\bf MSC2020 subject classifications}: Primary 60G42; 60G40; Secondary 60E15; 60F25.
\end{abstract}

\section{Introduction}

The study of martingales and their stopping times has been fundamental to
probability theory and its applications. In this paper, we examine two key
aspects of this relationship: the expectation of stopped martingales and
conditions needed to ensure that extended-valued stopping times are finite almost surely. The first of these results extends the classical optional stopping theorem by providing conditions under which the expected
value of a stopped martingale can be characterized by a finite real number, even
when the stopping time may be infinite with positive probability.

The second addresses a practical question in martingale theory: When can we guarantee that an extended-valued stopping time is finite almost surely? We present sufficient conditions that ensure this property and demonstrate their individual necessity through counterexamples. Our analysis reveals connections between the ``wobbling'' behavior of martingales and their stopping properties.

We also introduce three families of sums of i.i.d.~mean-zero random variables with distinct characteristics, demonstrating how our theoretical results apply in specific contexts. These families provide concrete examples in which the expected overshoot remains constant within each family, leading to explicit calculations of the aforementioned $L$.

\section{Historical Context}

The study of stopped martingales has a rich history dating back to Doob's fundamental work on martingale theory in the 1940s and 1950s. A foundational contribution to that development was Wald's 1944 computation of the expected value of a sequential probability ratio test involving a sum of i.i.d.\ random variables having a finite moment generating function. The number of summands was also random, determined by the number required for the partial sum $S_n$ to cross an upper boundary $A>0$ or a lower boundary $B<0$. The result came to be known as Wald's equation [13].

Soon after, Blackwell (1946) [1] extended and generalized Wald's lemma to the expectation of a sum of i.i.d.\ finite-mean random variables run up to a stopping time of finite expectation. Significant advances were made by Chow, Robbins, and Teicher (1965) [2] in their seminal paper on moments of randomly stopped sums, which established foundational results for the behavior of stopped sequences of martingale sums based on Wald's first and second equations. Later, Klass (1988) [10] and Klass (1990) [11] provided a best possible improvement of Wald's equation, obtaining definitive results for sums of independent, not necessarily identically distributed, mean-zero random variables.

Reinterpreting what was proved, let $X_1, X_2,\dots$ be a sequence of independent, not necessarily identically distributed, mean-zero random variables and let $T$ be a stopping time with respect to the martingale $S_n$, where $S_n=X_1+X_2+\dots+X_n$. Let $T^*$ denote an independent copy of $T$, also independent of all the $X_k$ for $k \geq 1$. Suppose $E|S_{T^*}|$ is finite. Order-of-magnitude approximations due to Klass (1981) [9] for fixed-$n$ random sums of independent random variables can be employed to determine whether $E|S_{T^*}|$ is finite. Then $E\max_{1 \leq n \leq T}|S_T|$ has the same order of magnitude. When this expectation is finite, the sequence $\{S_{T\wedge n}\}$ is uniformly integrable; therefore, $ES_T=0$. The theorem gives an explicit characterization of how long a stopping time $T$ can linger yet still force $E|S_{T'}|$ to be finite for all stopping times $T'$ such that $P(T' \geq n)\leq P(T \geq n)$ for all integers $n \geq 1$, while maintaining $ES_{T'}=0$. Victor de la Pe\~na would call $T^*$ a decoupled variable, characterizing this approach as based on decoupling. As noted by Grimmett and Stirzaker (2001) [8], ``Wald's equations provide a systematic approach to computing expectations involving random sums, with applications ranging from sequential analysis to financial mathematics'' (p.~418).

Decoupling methods then reached higher-order statistics. Chow, de la Pe\~na, and Teicher (1993) [4] extended Wald's equation to multilinear forms of arbitrary order---sums of products of $k$ distinct $X$'s---showing that the stopped form retains mean zero under a moment condition on the stopping time, and using this to bound the moments of the associated first-crossing times.\ de la Pe\~na and Lai (1997a) [5] treated the order-two case for general kernels, covering both regular and denormalized $U$-statistics and adding an asymptotic-bias analysis under random stopping.\ de la Pe\~na and Lai (1997b) [6] then unified and extended both, obtaining Wald's equation $EU_T=0$ for general degenerate kernels of arbitrary order, together with moment-convergence and asymptotic-expansion results. All of this work followed a single route---decoupling inequalities combined with H\"older's inequality---and de la Pe\~na and Zamfirescu (2002) [7] carried that route from $U$-statistics to general randomly stopped martingales, adding a domination technique to obtain a Wald identity in the martingale setting.

Let $S_{n,2}=\sum_{1\leq j\neq k\leq n} X_jX_k$. Again let $T$ denote a stopping time with respect to $S_{n,2}$. Assuming the $X_k$ were i.i.d.\ with mean zero, authors had shown that under various conditions, generally involving the finiteness of higher moments of the $X_k$ and/or $T$, $ES_{T,2}=0$ held. This appeared to be the universal theme. However, Klass (2009) [12] found a mean-zero $X$-distribution and a stopping time $T$ such that, if $T^*$ is an independent copy of $T$ and is also independent of the $X_k$, $E|S_{T^*\!,2}|$ is finite and yet $E|S_{T,2}|$ is infinite, depriving us of the ability to conclude that $ES_{T,2}$ is zero, or even that it has a well-defined Lebesgue integral. Since then, for various symmetric functions of independent mean-zero $X_1, X_2,\dots$, the ``Unsolved Problem'' has been: ``How slowly can a stopping time $T$ stop and yet retain the property that the expected value of its stopped sum has mean zero for all stopping times $T'$ such that $P(T'\geq n) \leq P(T\geq n)$?''

The work presented in this paper builds upon this historical foundation while extending the theory to address cases where traditional assumptions may not hold, particularly when stopping times can be infinite with positive probability. The traditional conclusion that $EM_TI(T<\infty)=EM_1$ need not hold.

\section{Main Results and Examples}

\subsection{Characterizing when \texorpdfstring{$EZ_T I(T<\infty)$}{EZT I(T<infinity)} is finite}

The following theorem can be used to determine, in the simplest possible way, whether the absolute value of a martingale run up to a possibly extended-valued stopping time $T$ has finite expectation on the set where $T$ is finite.

\begin{theorem}
Let $Z_n$ be any nonnegative sequence of ${\mathcal F}_n$-measurable random variables having finite expectation. Suppose that ${\mathcal F}_n$ is a nondecreasing sequence of sigma fields and $T$ is an associated extended-valued stopping time. Then $EZ_TI(T<\infty)<\infty$ if and only if $\lim EZ_TI(n \leq T < \infty)=0$.
\end{theorem}

{\bf Proof.}
Let $A=EZ_TI(T<\infty)$. First suppose $A<\infty$. Notice that $Z_TI(T\leq n)\to Z_TI(T<\infty)$ a.s. By monotone convergence, $\lim_{n\to\infty} EZ_TI(T\leq n)=A$. Hence $\lim_{n\to\infty} EZ_TI(n<T<\infty)=0$. Conversely, suppose $\lim EZ_TI(n\! <\!T\! <\!\infty)\! =\! 0$. Then there exists $n'$ such that
$EZ_TI(n'\! <\!T\! <\!\infty)\! \leq\! 1$. Therefore
\[
EZ_TI(T<\infty)\leq \sum^{n'}_{k=1} EZ_kI(T\!=\!k)+EZ_TI(n'\! < T \!<\!\infty)
\leq \sum^{n'}_{k=1} EZ_{k}+ 1<\infty.
\]
\hfill $\square$

\medskip\noindent
The utility of this theorem derives from the fact that whenever one needs to show that $EZ_TI(T<\infty)$ is finite, Theorem~3.1 characterizes the minimal fact that must be established. This is especially useful when $\lim EZ_{T\wedge n}$ is infinite.

\subsection{A conservation theorem obeyed by integrable stopped martingales}

\begin{theorem}
Let $\{M_n\}$ be a mean-zero martingale sequence and let $T$ be an extended-valued stopping time with respect to its filtration $\{{\mathcal F}_n\}$. Suppose $E|M_T|I(T<\infty)$ is finite. Then there is a finite real number $L$ such that $EM_T I(T<\infty)=L$ and $\lim_{n\rightarrow\infty} EM_nI(T>n)$ exists and equals $-L$. In particular, neither $L$ nor $P(T=\infty)$ need be 0.
\end{theorem}

This theorem extends to random sequences $\{M_n\}$ taking values in a Banach space.

\bigskip
\textbf{Proof.} By the dominated convergence theorem,
\begin{align*}
\lim_{n\rightarrow\infty}EM_T 1(T \leq n) = EM_T 1(T<\infty) \equiv L.
\end{align*}

Since $EM_{T \wedge n}=0$ and $M_{T \wedge n}$ splits as
\begin{align*}
M_{T \wedge n} = M_T 1(T \leq n) + M_n 1(T>n),
\end{align*}
we obtain
\begin{align*}
0 &= \lim_{n\rightarrow\infty}EM_T 1(T \leq n) + \lim_{n\rightarrow\infty}EM_n 1(T>n)\\
&= L + \lim_{n\rightarrow\infty}EM_n 1(T>n).
\end{align*}
Hence the limit of $EM_n 1(T>n)$ always exists and equals $-L$, proving the theorem. The same reasoning can be carried out in Banach spaces based on either norms or bounded linear functionals. \hfill $\square$

\bigskip
This result generalizes Theorem~4.1 of Klass (1988) [10] in two ways. First, the stopping time is allowed to be extended-valued, taking the value infinity with positive probability. Second, the new result only requires $E|M_T| I(T\!<\!\infty)\!<\!\infty$ instead of $\lim_{n\to \infty} E|M_{T\wedge n}|<\infty$, which was required by Theorem~2.3 of Chow, Robbins, and Siegmund (1971) [3]. The following example illustrates the extended applicability of Theorem~3.2.

\bigskip{\bf Example 1.} Let $M_n$ be a mean-zero martingale such that, for some finite positive $B$, $\limsup_{n\to\infty}|M_n|=\infty$ a.s.\ and $\liminf_{n\to\infty}|M_n|<B$ a.s. Let $T_0=0$. For all $k\geq 1$, let $T_{2k-1}$ be the first $n>T_{2k-2}$ such that $|M_n|>k^3$, and let $T_{2k}$ be the first $n>T_{2k-1}$ such that $|M_n|<B$. Then, for $k\geq 2$, let $T=T_{2k}$ with probability $k^{-2}$ and let $T=\infty$ otherwise. Then $E|M_T|I(T<\infty)<\infty$ and yet $\lim_{n\to\infty} E|M_{T\wedge n}|I(T<\infty)=\infty$. Letting $L=EM_TI(T<\infty)$, clearly $|L|<B P(T<\infty)$.

\begin{corollary}
For any $p>1$,
\begin{equation}
 |L|^{\frac{p}{p-1}} \leq\liminf_{n\rightarrow\infty} P(T>n)(E|M_n|^pI(T>n))^{\frac{1}{p-1}}.
\end{equation}
\end{corollary}

\textbf{Proof.} Let $U_n=M_nI(T>n)$ and $V_n=I(T>n)$. Notice that $V_n$ raised to any positive power is still $V_n$. By H\"older's inequality,
\begin{align*}
 |L|^{\frac{p}{p-1}} &\leq \liminf_{n\rightarrow\infty}\big( E|U_nV_n|\big)^{\frac{p}{p-1}}\\
& \leq \liminf_{n\rightarrow\infty} \big((E|U_n|^p)^{1/p} (EV_n)^{\frac{p-1}{p}}\big)^{\frac{p}{p-1}}\\
& = \liminf_{n\rightarrow\infty} EV_n(E|U_n|^p)^{\frac{1}{p-1}}.
\end{align*}
\hfill $\square$

\bigskip
This corollary is useful when $P(T<\infty)=1$. It also applies to Banach space random elements, provided the absolute values are replaced by the norm, or by bounded linear functionals. The corollary can be further extended, in utmost generality, by means of the Hausdorff--Young inequality.

\begin{corollary}
Using $p=2$, and assuming $EM_n^2=n$, for any extended-valued stopping time $T$ with $E|M_T|I(T<\infty)<\infty$,
$\liminf_{n\rightarrow\infty} nP(T>n) \geq L^2$. Recall that $L=EM_TI(T<\infty)$.
\end{corollary}

As with Corollary~3.2, this corollary is useful when $P(T<\infty)=1$.

\subsection{Three families of martingales with computable expected overshoot}

Here are three families of random variables that can be used as summands of three kinds of martingales with associated computable quantities.

\begin{itemize}
\item $C_1$ is the family of mean-zero integer-valued random variables whose maximum value is $+1$. For any random variable $X$ in $C_1$, $E(X|X>0)=1$.

\item For each $0<p<1$, $C(2,p)$ is the family of integer-valued mean-zero random variables $X$ whose positive integer values have probabilities of the form $P(X=k)=c(1-p)p^k$ for all integers $k>0$, for a suitable positive constant $c$ with $pc<1$. Notice that for each nonnegative integer $k$, $E(X|X>k)=k+\frac{p}{1-p}$\,.

\item For each $a\! >\! 0$, $C(3,a)$ is the family of mean-zero random variables such that for some $0<c<1$ and all $y>0$, $P(X\! >\! y)=c e^{-ay}$. For such $X$ and any $y>0$, $E(X|X>y)=y+\frac 1a$\,.
\end{itemize}

Associated with these collections of mean-zero random variables are three families of mean-zero martingales $F_1$, $F(2,p)$, and $F(3,a)$ whose increments given the past lie in $C_1$, $C(2,p)$, and $C(3,a)$, respectively. Let $h$ be any nonnegative real number, with $[h]$ denoting the largest integer $k\leq h$. Define the stopping time $T_h$ to be the first $n>0$ such that $M_n>h$, or $\infty$ if no such $n$ exists. Within each such family, the expected overshoot of any level is constant, so in $F_1$, $F(2,p)$, and $F(3,a)$, the value $E[M_{T_h} 1(T_h<\infty)]$ equals $([h]\!+\!1)P(T_h<\infty)$, $([h]\!+\!\frac {p}{1\!-\!p})P(T_h<\infty)$, or $(h\!+\!\frac 1a)P(T_h<\infty)$, respectively. Moreover, by the upcoming Theorem~4.1, it will follow that $P(T_h<\infty)=1$ for all three martingale families.

\subsection{Upper bounding the expected overshoot of exceedance levels}

The special character of these three martingale families enabled the calculation of the expected overshoot of any nonnegative exceedance level $r$ for the corresponding stopping time. We can abstract these conditions to obtain an explicit upper bound in a more general martingale framework.

\begin{theorem}
Take any $0<B<\infty$ and any nonnegative real number $r$. Let $\{M_n, {\mathcal F}_n\}$ be any martingale sequence.

{\rm{(i)}} Define
\[
T_r=
\begin{cases}
  \inf\{n \geq 1 : M_n > r\},\\
 \infty, \quad {\mbox{if no such finite $n$ exists.}}
\end{cases}
\]

{\rm{(ii)}} For all $n\geq 1$, suppose $E(M_n|T_r\!=\! n)\leq r+B$. Then
\begin{equation*}
E\{M_{T_r} I(T_r < \infty)\} \leq (r\!+\!B)P(T_r < \infty).
\end{equation*}
Equality holds if (ii) holds with equality for all $n$.
\end{theorem}

\section{On when a martingale has \texorpdfstring{$\limsup +\infty$}{limsup +infinity} almost surely}

\begin{theorem}
Let $\{M_n,{\mathcal F}_n\}$ be a mean-zero martingale. Take any nonnegative real number $r$. Define $T_r$ as the first $n\geq 1$ such that $M_n>r$ and set $T_r=\infty$ if no such $n$ exists. Suppose
\begin{enumerate}
\item[{\rm (i)}]
$\lim_{\varepsilon\downarrow 0} \limsup_{k\rightarrow\infty} P(\sup_{n\geq k}|M_n - M_k| > \varepsilon) = 1$.
\item[{\rm (ii)}] $\exists\, B < \infty$ such that for every $n \geq 1$,
$E(M_n|T_r = n)\leq r\!+\!B$ almost surely.
\end{enumerate}
Then $P(T_r<\infty)=1$ and $EM_{T_r} \leq r\!+\!B$. Thus Theorem~3.2 holds for $T=T_r$. Moreover, we may further conclude that $\limsup_{n\to \infty} M_n=\infty$ almost surely. Condition (i) ensures that the martingale wobbles arbitrarily far out in the sequence. Condition (ii) bounds the expected overshoot of $M_{T_r}$ beyond $r$.
\end{theorem}

\textbf{Proof.} Fix $r>0$. For $z>0$, let
\[
T_{r,z}=
\begin{cases}
\min\{n > 1 : M_n > r \ \text{or} \ M_n < -z\},\\
\infty, \quad \text{if no such } n \text{ exists.}
\end{cases}
\]

Write $M_n=X_1+\cdots+X_n$, set $Y_{j,z}=\max\{X_j,-r\!-\!z\}$, and let $M_{n,z}=\sum_{j=1}^n Y_{j,z}$. For each $z>0$, $M_{T_{r,z}\wedge n,z}$ is an $L^1$-bounded submartingale, so it converges almost surely. Hence the wobbles of the submartingale tend to 0.

Fix any sufficiently large $z>0$. We show that $P(T_{r,z}=\infty)=0$. To obtain a contradiction, suppose there exists $\delta>0$ such that $P(T_{r,z}=\infty)=2\delta$. By condition (i), there exist $0<\varepsilon<\min\{\delta, z/2\}$ and $k_0$ such that, for infinitely many $k\geq k_0$,
\begin{equation}
P\left(\sup_{n\geq k}|M_n - M_k| > \varepsilon\right) > 1 - \varepsilon.
\end{equation}

Since $M_{T_{r,z}\wedge n,z}$ converges almost surely, there exist $0<\varepsilon_1<\varepsilon$ and $k_1\geq k_0$ such that, for all $k\geq k_1$,
\begin{equation}
P\left(\sup_{n\geq k}|M_{T_{r,z}\wedge n,z} - M_{T_{r,z}\wedge k,z}| > \varepsilon_1\right) < \varepsilon.
\end{equation}

Hence, for all $k\geq k_1$,
\begin{align*}
2\delta\! -\! \varepsilon &< P\left(\{T_{r,z}=\infty\} \cap \{\sup_{n\geq k}|M_{T_{r,z}\wedge n,z} - M_{T_{r,z}\wedge k,z}| \leq \varepsilon_1\}\right)\\
&= P\left(\{T_{r,z}=\infty\} \cap \{\sup_{n\geq k}|M_{n,z} - M_{k,z}| \leq \varepsilon_1\}\right)\\
&= P\left(\{T_{r,z}=\infty\} \cap \{\sup_{n\geq k}|M_n - M_k| \leq \varepsilon_1\}\right)\\
&\leq P\left(\sup_{n\geq k}|M_n - M_k| \leq \varepsilon_1\right)\\
&\leq \varepsilon, \quad {\mbox{by inequality (2)}},
\end{align*}
so $2\delta<2\varepsilon$, a contradiction. Thus $P(T_{r,z}<\infty)=1$.

Let $P_{n,z}=P(M_{T_{r,z}\wedge n}>r)$ and $Q_{n,z}=P(M_{T_{r,z}\wedge n}<-z)$.

Since $P(T_{r,z}<\infty)=1$, $P_{n,z}+Q_{n,z}\to 1$ as $n\to\infty$. Moreover, $P_{n,z}$ is nondecreasing in $n$ with $P_{\infty,z}=P(M_{T_{r,z}}>r)=\lim_{n\to\infty}P_{n,z}$. Because $M_{T_{r,z}\wedge n}I(T_{r,z}>n)\leq r$,
\begin{equation*}
(r\!+\!B)P_{n,z}+rP(T_{r,z}>n)-zQ_{n,z}\geq EM_{T_{r,z}\wedge n,z} \geq EM_{T_{r,z}\wedge n}=0.
\end{equation*}
Sending $n\to\infty$,
\begin{equation*}
(r\!+\!B)P_{\infty,z} - z(1 - P_{\infty,z}) \geq 0.
\end{equation*}
Hence $P_{\infty,z}\geq \frac{z}{r+B+z}$. Observe that $P(T_r<\infty)\geq P_{n,z}$ for all $n\geq 1$ and $z>0$. Therefore $P(T_r<\infty)=1$. \hfill $\square$

\bigskip\noindent{\sc Comment.}
Theorem~4.1 can be extended to submartingales having nonnegative expectation.

\medskip
We now demonstrate that the conclusion of Theorem~4.1 can hold despite the fact that condition (ii) is violated.

\bigskip {\bf Example 3.}
Let $X_n$ be independent mean-zero two-valued random variables, with the two values being $2^{n-1}+n-1$ and $-2^{n-1}$. Notice that $P(X_n>0)$ converges to $\frac 12$. Let $T$ be the first $n$ such that $S_n=X_1+X_2+\dots+X_n>0$. Equivalently, $T$ is the first $n$ such that $X_n>0$. Since $P(T=n)$ decreases to zero exponentially fast, whereas $S_T|(T=n)=n$, it follows that $ES_T$ is positive and finite. Notice that condition (i) holds but condition (ii) fails. Nevertheless, by the Borel--Cantelli lemma, $X_n>0$ infinitely often, so $T$ is finite with probability 1.

\bigskip
The conditions of Theorem~4.1 do not guarantee that $\liminf M_n=-\infty$, even with positive probability.

\bigskip
{\bf Example 4.} (Illustrating that Theorem~4.1 can hold even when condition (ii) is not satisfied.) Let $X_1,X_2,\dots$ be independent mean-zero two-point random variables such that $X_n=1$ or $-n^2$. Then $X_n=1$ from some point on a.s. Hence the martingale $M_n=X_1+X_2+\dots+X_n$ tends to $\infty$ a.s. Furthermore, if $r$ is negative and $T'_r$ is the first $n\geq 1$ such that $M_n<r$, if such $n$ exists, and $T'_r=\infty$ otherwise, then $P(T'_r=\infty)>0$.

\bigskip A definitive statement can be made concerning when $EM_{T_r}<\infty$ for $r>0$: namely, iff
$\sum^{\infty}_{n=1} P(T_r=n)(r+{\mbox{the expected value of the martingale's overshoot of }} r)<\infty$\,.

\section{Two counterexamples to Theorem~4.1}

Conditions (i) and (ii) are critical. If either condition, or both conditions, fail, then Theorem~4.1 may not hold.

\bigskip{\bf Counterexample 1.}
(When (i) fails.) Let $X_1,X_2,\dots$ be independent mean-zero random variables such that $P(X_j=0)=1-j^{-2}$ and $P(|X_j|=1)=j^{-2}$. Then let
\[
T=\begin{cases}
\inf\{n\geq 1: X_1+X_2+\dots+X_n>0\},\\
\infty, \quad {\mbox{if no such $n$ exists.}}
\end{cases}
\]
Then (i) fails and (ii) holds, while $P(T=\infty)>0$.

This counterexample illustrates a phenomenon also noted by Grimmett and Stirzaker (2001) [8] in their discussion of recurrence properties of random walks. The distinction between recurrence and transience can depend on subtle properties of the increment distributions, particularly their tail behavior.

\bigskip{\bf Counterexample 2.} (When (ii) fails.)
Again let $X_1, X_2, \ldots$ be independent mean-zero random variables, but this time such that $P(X_j=(2^j\!-\!1)/j)=2^{-j}$ and $P(X_j=-1/j)=1-2^{-j}$. Then define $T$ as in Counterexample 1. Clearly, (ii) fails, (i) holds, and $P(T=\infty)>0$.

\section{Polynomial martingales}

Thus far we have applied our various theorems to what may be called linear martingales. There are also what can be considered square martingales, as well as polynomial ones. Our results also apply to them if they are martingales. As an example, suppose $S_n=X_1+X_2+\dots+X_n$ is a square-integrable mean-zero martingale. Then let $M_n=S_n^2-\sum_{j=1}^n X_j^2$. This $M_n$ is what we mean by a polynomial martingale of degree 2. A variation on this is $M_n^*=S_n^2\!-\!n$, assuming the conditional variance of every $X_j$ is 1.

The results presented in this paper extend classical martingale theory in two important directions: characterizing the expected value of stopped martingales when stopping times can be infinite, and providing sufficient conditions for ensuring that an extended-valued stopping time is finite almost surely. Our analysis builds upon and complements previous work by Wald (1944) [13], Blackwell (1946) [1], Klass (1988) [10], Klass (1990) [11], de la Pe\~na (1997a) [5], de la Pe\~na (2002) [7], and Klass (2009) [12], offering a unified perspective on the behavior of stopped martingales.

The three families of martingales introduced above provide concrete examples where our theoretical results can be applied, demonstrating the practical relevance of our findings. These examples also illustrate how different distributional properties of the increments can lead to distinct behaviors in the stopped process.

\section{Conclusion}

The theorems presented here reveal a precise structure underlying the apparent complexity of stopped martingales.

Theorems~3.1 and 3.2 apply to any mean-zero martingale and any extended-valued stopping time $T$, not just first-crossing times. They establish that the condition $\mathbb{E}|M_T|\mathbf{1}(T<\infty)<\infty$ is necessary and sufficient for the existence of a finite real number $L$ characterizing the stopped martingale's behavior---necessary because $L=\mathbb{E}M_T\mathbf{1}(T<\infty)$ would otherwise be undefined in the Lebesgue sense, and sufficient as shown by the proof using the dominated convergence theorem. When this condition holds, a fundamental conservation law emerges: the expected value $L$ on $\{T<\infty\}$ and the limiting expected value on $\{T>n\}$ sum to zero. The quantity $L$ that characterizes the instances that stop in finite time is exactly balanced by the shortfall, the amount $-L$ accumulated by those instances that do not. This duality persists whether $L$ is positive, negative, or zero, and whether $\mathbb{P}(T=\infty)$ is positive or zero---extending the classical optional stopping theorem into territory where the traditional assumptions fail.

For the important special case of first-crossing times, Theorem~4.1 identifies sufficient conditions under which $T=\inf\{n : M_n>0\}$ is finite almost surely: sufficient wobbling (condition (i)) combined with bounded overshoots (condition (ii)) guarantees that the martingale eventually crosses zero. The counterexamples demonstrate that both conditions are individually necessary in the sense that dropping either one can cause the conclusion to fail, though finding necessary and sufficient conditions for almost-sure finiteness of first-crossing times remains an open question.

\bigskip\noindent{\bf Acknowledgments.} Michael J.~Klass is Professor Emeritus of Statistics and Mathematics at the University of California, Berkeley, and wishes to acknowledge continued support. Victor de~la Pe\~{n}a is a Professor of Statistics at Columbia University, and acknowledges financial support from Google and DeepMind Project \#GTOO9019. We thank Fangyuan Lin for helpful comments following his careful reading of the paper.

\section*{References}
\begin{description}

\item{[1]} David Blackwell. On an extension of Wald's lemma.
 \emph{Ann. Math. Statist.} {\bf 17}, 310--315, 1946.

\item{[2]} Yuan~S. Chow, Herbert Robbins, and Henry Teicher.
 Moments of randomly stopped sums.
 \emph{Annals of Math. Stat.} {\bf 36}, 789--799, June 1965.

\item{[3]} Y.~S. Chow, H.~Robbins, and D.~Siegmund.
\emph{Great Expectations: The Theory of Optimal Stopping}, xii+139 pp.,
Boston, Houghton Mifflin, 1971.

\item{[4]} Y.~S. Chow, V.~H. de la Pe\~na, and H.~Teicher. Wald's equation for a class of denormalized $U$-statistics.
\emph{Ann.\ Probab.} \textbf{21}, 1151--1158, 1993.

\item{[5]} V.~H. de la Pe\~na and Tze Leung Lai (1997a).
Wald's equation and asymptotic bias of randomly stopped $U$-statistics.
\emph{Proc.\ Amer.\ Math.\ Soc.} \textbf{125}, 917--925, 1997.

\item{[6]} V.~H. de la Pe\~na and T.~L. Lai (1997b).
Moments of randomly stopped $U$-statistics. \emph{Ann.\ Probab.} {\bf 25} (4),
2055--2081, October 1997. \url{https://doi.org/10.1214/aop/1023481120}

\item{[7]} Victor H. de~la Pe\~{n}a and I.-M. Zamfirescu.
 Decoupling and domination inequalities with application to Wald's identity for martingales.
 \emph{Statist. Probab. Lett.} {\bf 57}, 157--170, 2002.

\item{[8]} G.~R. Grimmett and D.~R. Stirzaker.
  \emph{Probability and Random Processes}.
 Oxford University Press, Oxford, 3rd edition, 2001.

\item{[9]} Michael J. Klass.
 A method of approximating expectations of functions of sums of independent random variables.
\emph{Ann. Probab.} {\bf 9} (4), 413--428, 1981.

\item{[10]} Michael J. Klass.
 A best possible improvement of {W}ald's equation.
 \emph{Ann. Probab.} {\bf 16}, 840--853, 1988.

\item{[11]} Michael J. Klass.
 Uniform lower bounds for randomly stopped Banach space valued random sums.
 \emph{Annals of Probability}, {\bf 18}, 790--809, April 1990.

\item{[12]} Michael J. Klass.
 A denormalized U-statistic that cannot be decoupled from some associated stopping time.
 \emph{Statistics \& Probability Letters}, {\bf 79} (13), 1509--1511, 2009.

\item{[13]} Abraham Wald. On cumulative sums of random variables.
 \emph{Ann. Math. Statist.} {\bf 15}, 283--296, 1944.

\end{description}
\end{document}